\documentstyle[12pt]{article}
\begin{document}
\title{\bf Structures and Representations of Generalized Path Algebras }

\author{ {\em Shouchuan Zhang $^{a,~b}$ and Yao-Zhong Zhang $^b$ } \\
$a$. Department of Mathematics, Hunan University\\ Changsha 410082, P.R. China \\
$b$. Department of Mathematics, University of Queensland\\ Brisbane 4072, Australia }

\date{}

\newtheorem{Theorem}{\quad Theorem}[section]
\newtheorem{Proposition}[Theorem]{\quad Proposition}
\newtheorem{Definition}[Theorem]{\quad Definition}
\newtheorem{Corollary}[Theorem]{\quad Corollary}
\newtheorem{Lemma}[Theorem]{\quad Lemma}
\newtheorem{Example}[Theorem]{\quad Example}

\maketitle
\addtocounter{section}{-1}

\begin {abstract}
It is shown that an algebra $\Lambda $ can be lifted with nilpotent Jacobson radical  $r = r(\Lambda )$ and
  has  a generalized matrix  unit
$\{e_{ii}\}_I$ with each $\bar e_{ii} $ in the center of $\bar \Lambda = \Lambda /r$ iff $\Lambda $ is
isomorphic to a generalized path algebra with weak relations.
 Representations of the
generalized path algebras are given. As a corollary, $\Lambda $ is
a finite algebra with non-zero unity element over perfect field
$k$ (e.g. a field with characteristic zero or a finite field ) iff
$\Lambda $ is isomorphic to a generalized path algebra  $k (D,
\Omega , \rho )$ of finite directed graph with weak relations and
$dim {\ } \Omega < \infty $; $\Lambda $ is a generalized
elementary algebra which can be lifted with nilpotent Jacobson
radical and has a complete set of pairwise orthogonal idempotents
iff $\Lambda $ is isomorphic to a path algebra with relations.
\end {abstract}

\section {Introduction}

It is  well known that  every elementary algebra is isomorphic to
a path algebra of a finite directed graph with relations (see
\cite {ARS95}). In fact, every path algebra of a finite directed
graph with relations is also an elementary algebra. The results
are very useful because all representations of path algebras can
be obtained easily. In \cite {CL00} F.U. Coelho and S.X. Liu
introduced the concept of generalized path algebras to study other
algebras.

The aim of this paper is to give the structures and
representations of generalized path algebras  with weak relations.
We study   generalized path algebras by using generalized matrix
algebras introduced in \cite {Zh93}. In fact, every generalized
path algebra is a generalized matrix algebra. In section 1, we
study the structure of generalized matrix rings. We find the
relations among the decomposition of a ring, the complete set of
pairwise orthogonal idempotents (possibly infinite many) and
generalized matrix ring. This generalizes the theory about
decomposition of rings. In section 2, we study the representations
of the generalized path algebras. In section 3, we characterize
the generalized path algebras with weak relations by  algebras
which can be lifted with nilpotent Jacobson radical.

We say that an algebra $\Lambda $ can be lifted, if there exists a subalgebra
$A$ of $\Lambda $ such that $\Lambda = A \oplus r (\Lambda )$.
By the famous Wedderburn-Malcev Theorem (see \cite [Theorem 11.6 and Corollary
11.6] {Pi82}), for every finite dimensional
algebra $\Lambda $ over field $k$ with char $k=0 $, $\Lambda $ can be lifted
and $r(\Lambda )$ is nilpotent. We shall see, in section 3, that every
generalized path algebra with weak relations can be also
lifted and its Jacobson radical is nilpotent. In that section we show
that the converse also holds. That is,
it is shown that an algebra $\Lambda $ is isomorphic to a generalized path
algebra with weak relations iff $\Lambda $ can be lifted
with nilpotent Jacobson radical $r(\Lambda )$ and
has a complete set $\{e_{ii}\}_I$
of pairwise orthogonal idempotents with each $\bar e_{ii} $ in the center of $\bar \Lambda = \Lambda /r$.
As a corollary, $\Lambda $ is a finite algebra with non-zero unity element over field $k$ iff
$\Lambda $ is isomorphic to a generalized path
algebra  $k(D, \Omega , \rho )$ of finite directed graph  with weak relations and the dimension of $\Omega $ is finite; $\Lambda $ is a
generalized elementary algebra which can be lifted with nilpotent Jacobson
radical  iff $\Lambda $ is isomorphic to a path
algebra with relations.

\vskip.2in
\noindent ${\bf Preliminaries }$
\vskip.1in
Let $k$ be a field.
We first recall the concepts of $\Gamma _I$-systems, generalized matrix rings
(algebras ) and generalized path algebras.
Let $I$ be a non-empty set. If for any $ i, j, l, s \in I, A_{ij}$
is an additive group and there exists a map
$\mu_{ijl}$ from $A_{ij}\times A_{jl}$ to $A_{il}$ (written
$\mu _{ijl} (x, y)=xy)$ such that the following conditions hold:

(i) $(x +y)z = xz +yz, \ \ \ w(x+y)= wx +wy;$

(ii) $w(xz)=(wx)z$,

\noindent for any $x, y \in A_{ij}$, $z\in A_{jl}, w\in A_{li}$, then the
set $\{A_{ij }\mid i,j \in I\}$ is a $\Gamma _I$ -system with index $I$.

Let $A$ be the external direct sum of $\{ A_{ij} \mid i, j\in I \}$. We define
the multiplication in $A$ as $$xy = \{ \sum _k x_{ik}y_{kj} \}$$
for any $x=\{x_{ij}\}, y=\{y_{ij}\}\in A$.
It is easy to check that $A$ is a ring (possibly without the unity element ).
We call $A$ a generalized matrix ring, or a gm ring in short,
written as $A=\sum \{A_{ij} \mid i, j\in I\}.$ For any non-empty subset $S $ of
$A$ and $i, j \in I$, set $S_{ij} = \{ a\in A_{ij}
\mid \hbox { there exists } x \in S \hbox { such that } x_{ij} = a \}$.
If $B$ is an ideal of $A$ and $B = \sum \{B_{ij} \mid i, j \in I\}$, then $B$ is called a gm ideal. If for any $i, j \in I$,
there exists $ 0 \not= e_{ii} \in A_{ii}$ such that $x_{ij}e_{jj} =
e_{ii}x _{ij } = x_{ij}$ for any $x_{ij} \in A_{ij}$,
then the set $\{e_{ii} \mid i \in I\}$ is called a generalized matrix unit of
$\Gamma _ I$-system $\{A_{ij} \mid i, j \in I\}$, or a generalized matrix  unit of  gm ring $ A = \sum \{A_{ij}
\mid i, j \in I\}$, or a gm unit in short. It is easy to show that
if $A$ has a gm unit $\{e_{ii} \mid i \in I\}$, then every ideal $B$ of $A$ is
a gm ideal. Indeed, for any
$x = \sum _{i, j \in I } x_{ij} \in B$ and $i_0, j_0 \in I$, since $e_{i_0i_0}
x e_{j_0 j_0} = x_{i_0j_0} \in B$,
we have $B _{i_0j_0} \subseteq B$. Furthermore, if $B$ is a gm ideal of $A$, then  $ \{A_{ij} / B_{ij} \mid i,
j\in I\}$ is a $\Gamma _I$-system and
$A/ B \cong \sum \{A_{ij} / B_{ij} \mid i, j\in I\}$ as  rings.

If for any $ i, j, l, s \in I, A_{ij}$
is a vector space over field $k$ and there exists a $k$-linear map
$\mu_{ijl}$ from $A_{ij}\otimes A_{jl}$ into $A_{il}$ (written
$\mu _{ijl} (x, y)=xy)$ such that $x(yz)=(xy)z$ for any $x\in
A_{ij}$, $y\in A_{jl}, z\in A_{ls}$, then the set $\{A_{ij } \mid i,j \in I \}$
is a $\Gamma _I$- system with index $I$ over field $k$. Similarly, we get an
algebra $A=\sum \{A_{ij} \mid i, j\in I\},$
called a generalized matrix algebra, or a gm algebra in short.

Assume that $D$ is a directed (or oriented) graph ($D$ is possibly an infinite
directed graph and also possibly not a simple graph) (or quiver ). Let $I=D_0$
denote the vertex set of $D$ and $D_1$ denote the set of arrows of
$D$. Let $\Omega $ be a generalized matrix algebra over field $k$ with
gm unit $\{e_{ii} \mid i \in I \}$, the Jacobson radical
$r(\Omega _{ii})$ of $\Omega _{ii}$ is zero and $\Omega _{ij} =0$ for
any $i \not= j \in I$.  The sequence
$x=a_{i_0} x_{i_0i_1} a_{i_1}x_{i_1i_2}a_{i_2} x_{i_2i_3} \cdots x_{i_{n-1}i_
{n}} a_{i_n}$ is called a generalized path (or $\Omega $-path) from $i_0$ to
$i_n$ via arrows $x_{i_0i_1}, x_{i_1i_2},x_{i_2i_3}, \cdots ,
x_{i_{n-1}i_{n}}$,
where $0\not= a_{i_p} \in \Omega _ {i_p i_p} $ for $p = 0, 1, 2, \cdots, n $. In this case, $n$ is called the length of $x$,
written $l(x).$
For two $\Omega $-paths
$x=a_{i_0} x_{i_0i_1}a_{i_1}x_{i_1i_2}a_{i_2}x_{i_2i_3} \cdots
x_{i_{n-1}i_{n}} a_{i_n}$ and
$y=b_{j_0}y_{j_0j_1} b_{j_1}y_{j_1j_2}b_{j_2}y_{j_2j_3} \cdots
y_{j_{m-1}j_{m}}b_{j_m}$ of $D$ with $i_n=j_0$, we define the
multiplication of $x$ and $y$ as
$$
xy = a_{i_0}x_{i_0i_1}a_{i_1}x_{i_1i_2}a_
{i_2} x_{i_2i_3} \cdots x_{i_{n-1}i_{n}}(a_{i_n}b_{j_0})y_{j_0j_1}
y_{j_1j_2} b_{j_1}y_{j_2j_3} \cdots y_{j_{m-1}j_{m}}b_{j_m}.  \ \ \ \ \ \ \ (*)
$$

\noindent For any $i, j \in I,$ let $A_{ij}'$ denote the vector
space over field $k$ with basis being all $\Omega $-paths  from
$i$ to $j$ with length $>0$. $B_{ij}$ is the sub-space spanned by
all elements of forms:
 $$ a_{i_0} x_{i_0i_1}a_{i_1}x_{i_1i_2}a_{i_2}\cdots x_{i_{s-1}i_s} (a_{i_s}^{(1)} + a_{i_s}^{(2)} + \cdots  +
 a_{i_s}^{(m)})x_{i_s i_{s+1}}
 \cdots
x_{i_{n-1}i_{n}} a_{i_n} $$
$$ - \sum _{l = 1}^ma_{i_0} x_{i_0i_1}a_{i_1}x_{i_1i_2}a_{i_2}x_{i_2i_3} \cdots x_{i_{s-1}i_s} a_{i_s}^{(l)} x_{i_s i_{s+1}} \cdots
x_{i_{n-1}i_{n}} a_{i_n},$$
where $i_0 =i, i_n = j,  a _{i_s}^{(l)} \in \Omega _{i_s i_s},  a _{i_p} \in \Omega _{i_p i_p}$,
 $x_{i_t i_{t+1}} $  is an arrow, $p = 0, 1, \cdots , n$, $t = 0, 1, \cdots , n-1$, $l = 0, 1, \cdots , m$,
$0 \leq s \leq n$,  $n$ and $m$ are
 natural
numbers. Let $A_{ij} = A_{ij}'/B_{ij}$ when $i\not= j$ and $A_{ii}
= (A_{ii}' +\Omega _{ii})/B_{ii}$, written $[\alpha ] = \alpha
+B_{ij}$ for any generalized path $\alpha $ from $i$ to $j$.
 We can get a $k$-linear map from $A_{ij}\otimes A_{jl}$ to $A_{il}$  induced by
$(*).$ We
write $a$ instead of  $[a]$ when $a\in \Omega$. In fact, $[\Omega _{ii}] \cong \Omega _{ii}$ as algebras for any $i\in I.$
 Notice that we write $e_{ii}x_{ij} = x_{ij}e_{jj} = x _{ij}$ for any arrow $x_{ij}$ from $i$ to $j$. It is clear that
$\{A_{ij} \mid i, j \in I\}$ is a $\Gamma _I$-system with gm unit $\{e_{ii}
\mid i \in I\}$.
The gm algebra $\sum \{A_{ij}\mid i, j \in I \}$ is called the generalized path
algebra, or $\Omega $-path algebra, written as $k (D, \Omega )$
(see, \cite [Chapter 3]{ARS95} and \cite {CL00}). Let $J$ denote the ideal
generated by all arrows in $D$ of $k(D, \Omega )$.
If $\rho $ is a non-empty subset of $k (D, \Omega )$ and the ideal
$(\rho )$ generated  by $\rho $ satisfies $J^t \subseteq (\rho ) \subseteq J^2
$, then $k(D, \Omega )/(\rho )$ is called
generalized path algebra with relations. If $J^t \subseteq (\rho ) \subseteq J$,
then $k(D, \Omega )/(\rho )$ is called
generalized path algebra with weak relations. If $\Omega _{ii} = k e_{ii}$ for
any $i\in I,$ then $k (D, \Omega )$ is called a path algebra, written as $kD.$
If $D_0$ and $D_1$ are finite sets, then $D$ is called a finite directed graph.

Let $r(\Lambda )$ denote the Jacobson radical of ring $\Lambda $. Let $\mid $$S$$ \mid $
denote the number of elements in set $S$. Let $\delta _{ij}$ denote the
Kronecker $\delta $-function. Rings and algebras are possible without unity elements.

\section {Decomposition of generalized matrix rings}

In this section, we study the structure of generalized matrix rings. We find
the relations among the decomposition of a
ring, the complete set of pairwise orthogonal idempotents (possible
infinite many) and generalized matrix rings. This
generalizes the theory of direct sum decomposition of rings in \cite {AF92}.

\begin {Definition} \label {1.1}
If $A$ is a ring and $\{e_{ii} \mid i\in I \} \subseteq A$ such that the
following conditions are satisfied
(i) $e_{ii}e_{jj}= \delta _{ij}e_{ii}$ for any $i, j \in I$; (ii) for any $x\in
A$, there exists a finite subset $F$ of $I$ such that
$(\sum _{i\in F} e_{ii})x = x(\sum _{i\in F} e_{ii}) = x$; (iii) $e_{ii} \not=0$ for any $i\in I$,  then $\{e_{ii} \mid
i\in I \}$ is called the complete set of pairwise orthogonal
idempotents of  $A$ with index $I$. Moreover, if each  $e_{ii}$ is
a primitive idempotent (i.e. it can not be
written as a sum of two non-zero orthogonal idempotents), then
$\{e_{ii} \mid i\in I \} $ is called a complete set of pairwise orthogonal
primitive  idempotents of $A$ with index $I$
\end {Definition}

${\bf Remark:}$ (i) Let  $\{e_{ii} \mid i\in I \} $ be a complete set of pairwise orthogonal
 idempotents of $A$. Assume that  $x\in A $ and finite subset $F\subseteq I$ such that $x= (\sum _{i\in F} e_{ii})x=
x(\sum _{i\in F} e_{ii})=x$. If  $F' $ is a finite subset  of $I$ and $F\subseteq F'$, then
$x= (\sum _{i\in F'} e_{ii})x= x(\sum _{i\in F'} e_{ii})=x$. Indeed,
\begin {eqnarray*}
(\sum _{i\in F'} e_{ii})x = (\sum _{i\in F'} e_{ii})((\sum _{i\in F}
    e_{ii})x)
= ((\sum _{i\in F'} e_{ii})(\sum _{i\in F} e_{ii}))x
=  (\sum _{i\in F} e_{ii})x =x.
\end {eqnarray*}
Similarly, $x(\sum _{i\in F'} e_{ii})=x$.

(ii) Let $I$ be a non-empty set and  $A$  a ring with additive sub-groups $A_{ij}$ for any $i, j\in I$. If
$A = \sum _{i, j \in I}A_{ij}$ as additive groups and
$A_{ij} A_{st} \subseteq \delta _{js} A_{it}$ for any $i, j, s, t \in I$, then $\{A_{ij}, \mid i, j \in I\}$ is a $\Gamma _I$
-system. Let  $A'$ denote the gm ring  $ \sum \{A_{ij} \mid i, j \in I\}$ of
$\Gamma _I$-system  $\{A_{ij}, \mid i, j \in I\}$. Moreover, if $A_{ii}$ has  a non-zero  unity element $e_{ii}$ for any
$i \in I$, then  $A$ is the  inner direct sum of  $\{A_{ij}, \mid i, j \in I\}$ as additive groups
and $A'$ is isomorphic to $A$ under canonical isomorphism $\phi $ by sending $\{x_{ij}\}$ to $\sum _{i, j \in I}x_{ij}$
 for any $\{x_{ij}\} \in A'$. In this case,   $A$ is called  the inner gm ring of $\Gamma _I$-system
$\{A_{ij}, \mid i, j \in I\}$, also written  $A = \sum \{A_{ij}, \mid i, j \in I\}$.  If we view each element in $A_{ij}$ as  one in $ \sum \{A_{ij} \mid i, j \in I\}$,
then every
gm ring can be viewed as an inner gm ring.  Similarly, every
inner gm ring can be viewed as a  gm ring.

\begin {Theorem}\label {1.2}
$A$ has a complete set $\{ e_{ii} \mid i\in I \}$ of pairwise orthogonal
idempotents  with index $I$ iff $A = \sum \{ A_{i,j}\mid i, j \in I \}$
is a gm ring with gm unit $\{ e_{ii} \mid i \in I \}$ and $A_{ij} = e_{ii}A e_
{jj}$ for any $i, j \in I$.
\end {Theorem}
{\bf Proof.} The sufficiency is obvious. We now prove the necessity. Assume that  $A$ has a complete set
$\{ e_{ii} \mid i\in I \}$ of pairwise
orthogonal idempotents with index $I$. Let $A_{ij} = e_{ii}Ae_{jj}$ for any
$i, j \in I$. It is easy to check $A_{ij}A_{st} \subseteq \delta _{js} A_{it}$ for any $i, j, s, t \in I.$ Thus
$A$ is an inner gm ring of  $ \{ A_{i,j}\mid i, j \in I \}$ with gm unit $\{ e_{ii} \mid i \in I \}$.
$\Box$

This theorem implies that an algebra $A$ has a complete set of pairwise
orthogonal idempotents iff $A$ is a gm ring with gm unit.

\begin {Proposition} \label {1.3} (i) If $A$ has the non-zero   unity element $u$ then $A$ has a
complete set $\{ e_{ii} \mid i\in I \}$ of pairwise orthogonal
idempotents  with finite index $I$ and $\sum _{i\in I} e_{ii} =u.$

(ii) If ring $A$ has  the non-zero unity element $u$ and a complete set $\{ e_{ii} \mid i\in I \}$
of pairwise orthogonal idempotents  with index $I$, then $I$ is a finite set and
$\sum _{i\in I} e_{ii} =u.$

(iii) If $A$ is a finite dimensional algebra over field $k$, then $A$ has the non-zero unity element
iff $A$ has gm unit.
\end {Proposition}

{\bf Proof.} (i) Let $I = \{1\}$ and $u = e_{11}$.

(ii) Since $A$ has a gm unit $\{e_{ii}\}_I$, by Theorem \ref {1.2}, $A = \sum
\{A_{ij} \mid i, j \in I\}$ is a gm ring
with gm unit $\{e_{ii}\}_I$ and $A_{ij} = e_{ii} A e_{jj}$ for any $i. j \in
I.$ Let $u= \sum _{i, j\in F} u_{ij}$ with finite subset
$F$ of $I$ and $u_{ij} \in A_{ij}$ for any $i, j \in F$. Since $u$ is the  unity element
of $A$, $A_{ij} = 0$ for any $i \not\in F$ or $j \not\in F$. Thus $F= I$ since $e_{ii} \not=0$ for any $i\in I.$
For any $s \in I$ and $x_{ss}\in A_{ss}$,
since $u x _{ss} = x_{ss} $ and $ x_{ss}u = x_{ss} $, we have
$u_{ss}x_{ss} =x_{ss}$ and $x_{ss}u_{ss}=x_{ss} $.
This implies $u_{ss} = e_{ss}$ for any $s\in F.$
Next we show $u_{ij}=0$ when $i \not= j.$ On the one hand, $u _{ii}u = u_{ii}$.
On the other hand, $u _{ii} u = \sum _{s\in I}u_{ii}u_{is}$. Consequently, $u_
{ij}=0$ for any $i\not=j.$

(iii) If $A$ has gm unit $\{e_{ii}\}_I$, then  $I$ is finite since $A$ is finite dimensional. It is clear that
$u = \sum _{i \in I}e_{ii}$ is the unity element of $A$. The converse follows from (i). $\Box$

\begin {Proposition} \label {1.4}
 If  $A$ is a left (or right ) artinian  or noetherian ring with gm unit $\{e_{ii}\}_I$,
 then $I$ is finite and $\sum _{i\in I} e_{ii}$ is  the unity element of $A$.
\end {Proposition}

{\bf Proof.}  By Theorem \ref {1.2}, $A = \sum \{A_{ij} \mid i, j \in I\}$ with $A_{ij} = e_{ii}Ae_{jj}$ for
 any $i, j\in I.$  If $I$ is infinite, then there exists an infinite  sequence $i_1, i_2, \cdots,  i_n, \cdots $ in $I$.
Let $A_1 = Ae_{i_1i_!}, $ $A_2 = A_1 + Ae _{i_2i_2}, \cdots , A_{n+1} = A_n + A_{i_{n+1}i_{n+1}}, \cdots$.
Obviously   $A_1 \subset A_2 \subset \cdots \subset A_n \subset
\cdots $ is an ascending chain of left ideals of $A$.
 Let $B_1 = \sum _{j \in I, j\not= i_1} Ae_{jj},
B_2 = \sum _{j \in I, j\not= i_1, i_2} Ae_{jj},  \cdots , B_{n+1} = \sum _{j \in I, j\not= i_1, i_2, \cdots , i_{n+1}}
 Ae_{jj}$ for any natural number $n.$ Obviously,
 $B_1 \supset B_2 \supset \cdots \supset B_n \supset \cdots $ is an descending chain of left ideals of $A$.
 We get a contradiction.
 Consequently, $I$ is finite.
$\Box$

Let ${\cal A} \Gamma _I$ denote the category of all $\Gamma _I$ -systems with gm
unit, the morphism of two objects from $\{ A_{ij} \mid i, j \in I\}$ with
gm unit $\{e_{ii}\}_I$ to $\{B_{ij} \mid i,j\in I\}$ with gm unit $\{e_{ii}'\}_I$
is a set $\{f _{ij}\}_I$, where $f_{ij}$ is an additive group homomorphism
from $A_{ij}$ to $B_{ij}$ with $f_{ij} (xy) = f_{is}(x)f_{sj}(y)$ and
$f_{ii}(e_{ii}) = e_{ii}'$ for any $i, j, s\in I, x\in A_{is}, y \in A_{sj}.$
Let ${\cal GM}_I$ denote the category of all generalized matrix algebras with
index $I$ and gm unit, the morphism between the two objects is gm
homomorphism. A gm homomorphism of two objects from
$A = \sum \{ A_{ij} \mid i, j \in I\}$ with gm unit $\{e_{ii}\}_I$ to $B= \sum
\{B_{ij} \mid i, j \in I\}$ with gm unit $\{e_{ii}'\}_I$ is a ring homomorphism
$f: A \rightarrow B $ such that $f (A_{ij})\subseteq B_{ij}$ and
$f(e_{ii}) =e_{ii}'$ for any $i, j \in I.$

\begin {Proposition} \label {1.5} ${\cal A} \Gamma _I$ and ${\cal GM}_I$ are
two equivalent categories.
\end {Proposition}

{\bf Proof.} Let $H: {\cal A} \Gamma _I \rightarrow {\cal GM}_I$ by
$H(\{A_{ij}\}_I) = \sum \{ A_{ij} \mid i, j \in I\},$
$H(\{f_{ij}\}_I) = \oplus _{i, j \in I} f_{ij} $ for any morphism
$\{ f_{ij} \}_I$ from $\{ A_{ij} \mid i, j \in I\}$ to
$ \{B_{ij} \mid i, j \in I\}$. Let $G: {\cal GM}_I\rightarrow {\cal A}\Gamma
_I$ by $G( \sum \{ A_{ij} \mid i, j \in I\}) = \{A_{ij}\}_I$ and
$G(f)= \{f_{ij}\}_I$ with $f_{ij} = f$$\mid _{A_{ij}}$ for any $i, j \in I$.
Obviously, $HG = id$ and $GH= id$. $\Box$

\section {Representations of generalized path algebras}

In this section, we study representations of the generalized path algebras.
%We assume that $A$ is a generalized matrix algebra with the gm unit
%$\{e_{ii}\}_I$. Obviously, every generalized path algebra has gm unit.

\begin {Definition} \label {2.1} Let $\{ A_{ij} \mid i, j \in I\}$ be an
$\Gamma_I$-system with gm unit $\{e_{ii}\}_I$.
For any $ i, j \in I, M_{i}$ is an additive group and there exists a map $
\phi _{ij}$ from $A_{ij} \times M_j$ to $M_i$ (written
$ \phi _{ij} (a, x) = ax$) such that the following conditions are satisfied:

(i) $a(x+y)=ax+ ay$ and $(a+b)x=ax+bx.$

(ii) $(ca)x=c(ax).$

(iii) $e_{jj}x=x$

\noindent For any $x, y\in M_{j}, a, b\in A_{ij}, c\in A_{si}$,
then $\{M_i \mid i \in I\}$ is called an $\{A_{ij}\}_I$- module system.
\end {Definition}

Let {\cal R}ep $\{ A_{ij}\}_I$
denote the category of $\{ A_{ij}\}_I$-module systems. The morphism of two objects $\{M_{i}\}_I$ and $\{N_{i}\}_I$ is
a  collection  $\{f_{i}\}
_I$  such that $f_i$ is an additive group homomorphism from $M_i$ to $N_i$ with
$f_i(a_{ij}x_{j}) = a_{ij}f_j(x_j)$ for any $a_{ij}\in A_{ij}, x_j \in M_j.$

An $A$-module is called a local unitary $A$-module if for any $x\in M$ there
exists $u \in A$ such that $ux=x$.

\begin {Lemma}\label {2.2} If $A$ is a gm ring with gm unit $\{e_{ii}\}_I$,
then $M$ is a local unitary $A$-module iff $M$ is an $A$-module with $AM =M$.
\end {Lemma}

{\bf Proof.} Assume $AM = M$. For any $x\in M,$ there exist $a^{(p)} \in A, x^
{p}\in M$ such that $x = \sum _{p=1} ^n a ^{(p)}x^{(p)}$. There exists
a finite subset $F$ of $I$ such that $a^{(p)} \in \sum _{i, j \in F} A_{ij}$ for $p =1, 2, \cdots , n$. Let
$u = \sum_{i \in F} e_{ii}.$
We have that $u x = u(\sum _{p = 1, 2, \cdots , n } a^{(p)} x^{(p)}) =
\sum _{p = 1, 2, \cdots , n } a^{(p)} x^{(p)} =x$. Therefore, $M$ is a local unitary $A$-module.
Conversely, it is clear that $AM =M$ when $M$ is a local unitary $A$-module. $\Box$

\begin {Lemma}\label {2.3} Let $A$ be a gm ring with gm unit $\{e_{ii}\}_I$.

(i) If $M$ is a local unitary $A$-module, then $\{M_{i} \mid i \in I \}$ is an
$\{ A_{ij} \}_I$-module system with $e_{ii} M= M_i $.

(ii) If $\{M_{i}\}_I$ is an $ \{ A_{ij}\}_I$-module system, then the external
direct sum $M$ of $\{M_{i}\}_I$ becomes a local unitary $A$-module
under module operation $a x =\{ \sum _{s \in I} a_{is}x_{s}\}_I$ for any $a= \{
a_{ij}\}_I\in A$, $x = \{x_{i}\}_I\in M$.
\end {Lemma}

{\bf Proof.} (i) If $M$ is a local unitary $A$-module. Set $e_{ii} M= M_i$ for
any $i \in I$. It is clear that $\{M_{i}\}_I$ is an $\{ A_{ij}\}_I$-module
system. Indeed, for any $x, y\in M_{j}, a, b\in A_{ij}$ and $ c\in A_{si}$,
we have that $a(x+y)=ax+ ay$, $(a+b)x=ax+bx$, $(ca)x=c(ax)$ and $e_{jj}x=x$.

(ii) It is clear. Indeed, for any $a= \{ a_{ij}\}_I, b= \{b_{ij}\}_I\in A$ and
$x = \{x_{i}\}_I\in M$, it is easy to check $(ab)x = a(bx)$. Since there
exists finite subset $F$ of $I$ such that $x= \sum _{i \in F} x_{i}$, we have
that $ (\sum _{i \in F} e_{ii})x =x$. Thus $M$ is a local unitary $A$-module.
$\Box$

Let ${\ }_A{\cal M}LU$ denote the category of local unitary $A$-modules. every
morphism of two objects $M$ and $N$ is a homomorphism of $A$- modules.

\begin {Theorem}\label {2.4} Let $A = \sum \{A_{ij} \mid i, j \in I\}$ be a gm
ring with gm unit.
Then {\cal R}ep $\{ A_{ij}\}_I$ and ${}_A {\cal M}LU$ are equivalent.
\end {Theorem}

{\bf Proof.} Let $H:$ {\cal R}ep $\{ A_{ij}\}_I \rightarrow$ ${}_A {\cal M}LU$
by $H( \{M_i\}_I) = \sum \{M_{i} \mid i \in I\},$
$H(\{f_{i}\}_I) = \oplus _{i \in I}f_{i} $ for any morphism $\{f_{i}\}_I$
between two objects $\{M_{i}\}_I$ and $\{N_{i}\}_I$.
Let $ G:$ ${}_A{\cal M}LU \rightarrow $ {\cal R}ep $\{ A_{ij}\}_I $ by
$G(M)= \{M_{i}\}_I$ with $M_i = e_{ii}M$ for any $i\in I$. $G(f) = \{f_i\}_I $
with $f_i = f$$\mid _{M_i}$ for any morphism $f $ between two objects $M$ and
$N$. It is clear $HG = id$ and $GH = id.$ $\Box$

If $A = \sum \{A_{ij} \mid i, j \in I\}$ is a gm algebra over field $k$ with gm
unit $\{e_{ii}\}_I,$ we can similarly define  $\{A_{ij}\}_I$-module systems as follows.

Let $\{ A_{ij} \mid i, j \in I\}$ be a $\Gamma _I$ -system over field $k$ with
gm unit $\{e_{ii}\}_I$. If for any $ i, j \in I, M_{i}$ is a vector space
and there exists $k$-linear map $ \phi _{ij}$ from $A_{ij} \otimes M_j$ to
$M_i$ (written $ \phi _{ij} (a, x) = ax$) such that the following conditions
are satisfied:

(i) $(ca)x=c(ax).$

(ii) $e_{jj}x=x,$

\noindent for any $x \in M_{j}, a \in A_{ij}, c\in A_{si}$,
then $\{M_i \mid i \in I\}$ is called an $\{A_{ij}\}_I$- module system.
We still use the two notations {\cal R}ep $\{ A_{ij}\}_I$ and
${}_A {\cal M}LU$ to denote the corresponding categories.

\begin {Theorem}\label {2.5} Let $A = \sum \{A_{ij} \mid i, j \in I\}$ be a gm
algebra with gm unit.
Then {\cal R}ep $\{ A_{ij}\}_I$ and ${}_A {\cal M}LU$ are equivalent.
\end {Theorem}

For a generalized path algebra $k (D, \Omega , \rho )$ with weak relations, let
$P = k (D, \Omega)$, $N = (\rho ) $ and $Q = P/ N.$
It is clear that the generalized path algebra $k (D, \Omega , \rho )$ with weak relations  is a gm
algebra, so its representation corresponds to
$\{Q_{ij}\}_I$-module system. That is, {\cal R}ep $\{ Q_{ij}\}_I$ and ${}_Q
{\cal M}LU$ are equivalent. However, we have  a simpler category.

A representation of $(D, \Omega )$  is a set $ (V, f) =: \{ V_i,
f_\alpha \mid V_i$ is an unitary $\Omega _{ii}$-module,
$f_\alpha : V_i\rightarrow V_j$ is a $k$-linear map, $ i, j \in I$,  $ \alpha $ is
an arrow from $j$ to $i \}$.
A morphism $h: (V, f) \rightarrow (V', f')$ between tow representations of $(D, \Omega )$
is the collection $\{h_i\}_I$ such that $h_i : V_i \rightarrow V_i'$ is a
$k$-linear map and $h_j f _\alpha  = f_\alpha 'h_i$ for any arrow
$\alpha : i \rightarrow j$ and $i, j \in I.$
Let {\cal R}ep $(D, \Omega )$ denote the category of representations of $(D,
\Omega )$.

\begin {Lemma}\label {2.6} Let $P = k (D, \Omega )$ and $Q = k(D, \Omega , \rho ).$

(i) If $(V, f)$ is an object in
{\cal R}ep $(D, \Omega )$, then $\{V_i\}_I$ is a $\{ P_{ij}\} _I$-module
system under operation $ \alpha \cdot v_{i_n} =
a_{i_0}\cdot f _{x_{i_0i_1}} ( a_{i_1} \cdot (f_{ x_{i_1i_2}} \cdots
f_{x_{i_{n-1}i_n}}(a_{i_n} \cdot v_{i_n}) ))$ for any $\Omega $-path
$\alpha = a_{i_0}x_{i_0i_1}a_{i_1} x_{i_1i_2} \cdots x_{i_{n-1}i_n} a_{i_n}$
from $i_0$ to $i_n$ and $v_{i_n} \in V_{i_n}$.

(ii) If $\{V_i \}_I$ is a $ \{P_{ij}\}_I $-module system, then $(V, f)$ is an
object in {\cal R}ep $(D, \Omega )$ under operation
$f_{x_{ij}} (v_j) = x_{ij}\cdot v_j $ for any arrow $x_{ij} \in P_{ij}$ and
$v_j \in V_j.$

\end {Lemma}

{\bf Proof.} (i) It is sufficient to show that
$$
(\alpha \beta ) \cdot v_{j_m} = \alpha \cdot (\beta \cdot x_{j_m}) \ \ \ \ \
\ \ \ \ \ \ \ \ \ \ (*)
$$
for two $\Omega $- paths
\ \ \ \ \ \ \ \ \ \ $\alpha =a_{i_0} x_{i_0i_1}a_{i_1}x_{i_1i_2}a_{i_2}x_
{i_2i_3} \cdots x_{i_{n-1}i_{n}}a_{i_n}$
\ \ \ \ \ \ \ \ and \\
$\beta =b_{j_0}y_{j_0j_1} b_{j_1}y_{j_1j_2}b_{j_2}y_{j_2j_3} \cdots
y_{j_{m-1}j_{m}}b_{j_m}$ of $D$ with $i_n=j_0$.

When $\alpha \beta \not=0,$ i.e. $a_{i_n} b_{j_0} \not=0,$ $\alpha \beta $ is
an $\Omega$-path. By definition, (*) holds.
When $\alpha \beta =0,$ i.e. $a_{i_n} b_{j_0} =0,$ $\alpha \beta $ is not an
$\Omega$-path. Obviously the left side of (*) =0.
\begin {eqnarray*}
\hbox { The right side of } (*) &=& \alpha \cdot
(b_{j_0}\cdot f _{y_{j_0j_1}} ( b_{j_1} \cdot f_{ y_{j_1j_2}}( \cdots f_{y_{y_
{m-1}j_m}}(b_{i_m} \cdot v_{i_m}) )))\\
&=& a_{i_0}\cdot f _{x_{i_0i_1}} ( a_{i_1} \cdot f_{ x_{i_1i_2}}( \cdots \\
&{\ }&f_{x_{i_{n-1}i_n}}((a_{i_n}b_{j_0})\cdot f _{y_{j_0j_1}} ( b_{j_1}
\cdot f_{ y_{j_1j_2}}( \cdots f_{y_{y_{m-1}j_m}}(b_{i_m} \cdot v_{i_m}) )))))\\
&=& 0.
\end {eqnarray*}
Consequently, (*) holds.

(ii) It is obvious. $\Box$

Combining Lemma \ref {2.6} and Theorem \ref {2.5}, we have
\begin {Theorem}\label {2.7}
{\cal R}ep $(D, \Omega )$ and ${}_{k(D , \Omega ) } {\cal M}LU$
are equivalent.
\end {Theorem}

For a representation $(V, f)$ in {\cal R}ep $(D, \Omega )$ and any element
$\sigma \in k(D, \Omega )$, by Lamma \ref {2.6} and Theorem \ref {2.5}, $(V, f)$  can be viewed as $k(D, \Omega )$-module,
so for any $\sigma \in k(D,\Omega )$,
we  write  $f_\sigma : V \rightarrow V$ by sending $x$ to $\sigma \cdot x$ for any $x\in V$.
Let {\cal R}ep $(D, \Omega , \rho )$ denote the full subcategory of {\cal R}ep
$(D, \Omega )$ whose objects are $(V, f)$ with $f _\sigma =0$ for each
$\sigma \in \rho .$

\begin {Lemma}\label {2.8} Let $P = k (D, \Omega )$ and $Q = k(D, \Omega , \rho ).$

(i) If $(V, f)$ is an object in {\cal R}ep $ (D, \Omega , \rho )$, then
$\{V_i\}_I$ is a $\{ Q_{ij}\} _I$-module system under operation induced by
operation of $\{P_{ij}\}$- module system in Lemma \ref {2.6}.

(ii) If $\{V_i \}_I$ is a $ \{Q_{ij}\}_I $-module system, then $(V, f)$ is an
object in {\cal R}ep $(D, \Omega , \rho )$ under operation
$f_{x_{ij}} (v_j) = x_{ij}\cdot v_j $ for any arrow $x_{ij} \in P_{ij}$ and
$v_j \in V_j.$

\end {Lemma}

\begin {Theorem}\label {2.9}

(i) {\cal R}ep $(D, \Omega ,\rho )$ and ${}_{k(D , \Omega , \rho )
} {\cal M}LU$ are equivalent.

(ii) If $D$ is finite (i.e. $I$ is finite and the number of arrows
between any two vertexes is finite ), then f.d.{\cal R}ep $(D,
\Omega, \rho )$ and f.d.${}_{k(D , \Omega , \rho )}{\cal M}LU$ are
equivalent. Here, f.d.{\cal R}ep $(D, \Omega , \rho )$ and
f.d.${}_{k(D , \Omega , \rho ) } {\cal M}LU$ denote the full
subcategories of finite dimensional objects in the corresponding
categories, respectively.

\end {Theorem}

\section {Generalized path algebras }

In this section, we characterize the generalized
path algebras with weak relations by some
algebras which can be lifted with nilpotent Jacobson radical.

If $V = U \oplus W $ as vector spaces and $x\in V $, then there exist $a\in U$ and $b\in W$ such that
$x = a + b.$ For convenience, we denote $a$ and $b$ by $x_U$ and $x_W,$ respectively.

\begin {Lemma} \label {3.1} Let $\Lambda $ be an algebra and $N$ an ideal of
$\Lambda $. Then the following conditions are equivalent:

(i) There exists a subalgebra $A$ of $\Lambda $ such that
$\Lambda = A \oplus N$ as vector spaces.

(ii) The canonical homomorphism $\pi : \Lambda \rightarrow \Lambda /N$ is
split in the category of algebras, i.e. there exists an algebra homomorphism
$\xi : \Lambda /N \rightarrow \Lambda $ such that $\pi \xi = id _{\Lambda /N}.$

\end {Lemma}

{\bf Proof.}
(i) $\Rightarrow $ (ii).
Define $\xi : \Lambda /N \rightarrow \Lambda $ by sending $\xi (x +N) = x_A$
for any $x = x_A +x_N \in \Lambda $ with  $x_A \in A, x_N \in N.$ It is clear that
$\xi $ is an algebra homomorphism and $\pi \xi = id.$

(ii) $\Rightarrow $ (i). Obviously $\Lambda = A \oplus N$ with $A = Im \xi.$
$\Box$

We say that an algebra $\Lambda $ can be lifted if $\Lambda = A
\oplus r(\Lambda )$  with subalgebra $A$ {\ }\footnote {This
concept was introduced by Fang Li. }.

\begin {Lemma} \label {3.2}

Let $\Lambda $ be an algebra, $N$ an ideal of $\Lambda $ and $A$ a subalgebra
of $\Lambda$. If $\Lambda = A \oplus N$, then
$\Lambda /B = (A+B)/ B \oplus (N+B)/B$ for any ideal $B$ of $\Lambda $ with
$B \subseteq A$ or $B \subseteq N$.

%(ii) If there exists an idempotemt $0\not= e\in \Lambda $ such that $\Lambda /
%N \cong ke$ as algebras, then $\Lambda = ke \oplus N$ can be lifted.

\end {Lemma}

{\bf Proof.}
For any $x = x_A + x_N \in \Lambda $ with $x_A\in A$ and $x_N \in N$,
$\bar x =x +B =(x_A +B) + (x_N +B) \in \Lambda /B$ with
$(x_A +B) \in (A +B)/B, (x_N +B)\in (N+B)/B$. This implies that $\Lambda /B=
(A +B)/B + (N +B)/B$. Assume $B \subseteq A$. then $(A/B)\cap ((N+B)/B)=0$
and $\Lambda /B = A/B\oplus (N+B)/B.$
Similarly, when $B \subseteq N$, $\Lambda /B = (A+B)/B \oplus N/B.$ $\Box$
\begin {Lemma} \label {3.3}

Let $\Lambda $ be an algebra, $N$ a nilpotent  ideal of $\Lambda $ and $A$ a subalgebra
of $\Lambda$. Assume  $\Lambda = A \oplus N$ as vector spaces. If  $\{e_{ii}\}_I$ is a complete set of
pairwise orthogonal
idempotents of $\Lambda $, then $\{e_{ii} \}_I \subseteq A. $
\end {Lemma}
{\bf Proof.} We first show that if $e$ is  idempotent in $\Lambda $ with $e = e_A + e_N$ and $e_A \in A, e_N \in N$,
then $e_A$ is idempotent.  Indeed, since $ee=e$ and $N$ is an ideal of $\Lambda$,
we have $e_A e_A + (e_A e_N + e_N e_N +e_Ne_A) = e_A + e_N,$
 which implies that $e_A e_A = e_A$.

Next we show that if $e$ and $f$ are pairwise orthogonal idempotents of $\Lambda $, then so are  $e_A$ and $f_A$.
Indeed, since $ef =0 $, i.e. $e_A f_A + (e_Af_N + e_Nf_A + e_N f_N) =0$, we have $e_A f_A =0.$ Similarly, $f_Ae_A=0.$

We now show that  each $e_{ii} \in A$
  by induction for $m$, where $N^m =0$.

When $m=1$, $N=0$. In this case, $(e_{ii})_A = e_{ii}\in A$ for any $i\in I.$

Assume now that the claim holds when $m \leq l$ and we show that the claim
also holds when $m= l+1.$
Let $\bar \Lambda = \Lambda /N^l.$ By Lemma \ref {3.2}, $\bar \Lambda = (A + N^l)/N^l \oplus N/N^l$. It is clear
$\{\bar e_{ii}\}_I$ is a complete set of pairwise orthogonal idempotents of $\Lambda /N^l.$ By the inductive assumption,
$\bar e_{ii} \in \bar A$,  i.e. $(e_{ii})_N \in N^l $ for any $i\in I.$

For any $x \in \Lambda $, there exists a finite subset $F$ of $I$ such that
\begin {eqnarray}
x = (\sum _{i\in F} e_{ii})x
\hbox { \ \ \ \ \ \ \ and  \ \ \ \ \ \ \ \ \ \ } x_A = (\sum _{i\in F} e_{ii})x_A. \label {e3.2.12}
\end {eqnarray}
By (\ref {e3.2.12}),
\begin {eqnarray}
0 = (\sum _{i\in F}  (e_{ii})_ N )  x_A  \hbox { \ \ \ and \ \ \ } x_A = (\sum _{i\in F} (e_{ii})_A)x_A   \label {e3.2.15}.
\end {eqnarray}
Since $(\sum _{i\in F}  (e_{ii})_ N )  x_N  \in N^{l+1} =0 $, $(\sum _{i\in F}  (e_{ii})_ N )  x_N  =0$.
  By (\ref {e3.2.12}) and (\ref {e3.2.15}),
\begin {eqnarray}
 x _N = (\sum _{i\in F}  (e_{ii})_ A )  x_N  \label {e3.2.16}.
\end {eqnarray}
Combining   (\ref {e3.2.15}) and ( \ref {e3.2.16}), we have that  $ x = (\sum _{i\in F}  (e_{ii})_ A )  x$.
Similarly, $ x = x(\sum _{i\in F}  (e_{ii})_ A )  $. Consequently,
 $\{ (e_{ii})_{ A}\}_I$ is a complete set of pairwise orthogonal idempotents of $\Lambda .$
Since $e_{ii}$ and $(e_{ii})_A $ are the unity element of $\Lambda _{ii}$, $e_{ii} = (e_{ii})_A\in A$ for any $i\in I.$
$\Box$

By Lemma \ref {3.3}, we have immediately:
\begin {Lemma} \label {3.4}
Let $\Lambda $ be an algebra with non-zero unity element  $u$, $N$ a nilpotent  ideal of $\Lambda $ and $A$ a subalgebra
of $\Lambda$. If   $\Lambda = A \oplus N$ as vector spaces, then $u \in A$.
\end {Lemma}

\begin {Lemma}\label {3.5} Let $A$ be a subalgebra of $\Lambda $ and
$\Lambda = A \oplus r $ with nilpotent Jacobson radical $r= r (\Lambda )$.
Let $ B = \{ r_u \mid u \in U \} \subseteq r $. If $\bar B = \{ \bar r_u \mid u
\in U \} $ generates $r/r^2$ as $\Lambda /r$-modules, then $A\cup B$
generates $\Lambda $ as algebras.
\end {Lemma}

{\bf Proof.}
Since $r$ nilpotent, there is  $m$ such that
$r^m =0.$  We use induction on $ m.$
It is obvious that  $r=0$ and $\Lambda =A$ when $m=1$.
When $m=2$, we have that $r^2 =0$ and $r= r/r^2$. Thus $\bar B =B$ generates
$r$ as $\Lambda /r$-modules. That is,
$r = \sum _{u\in U}\Lambda r_u= \sum _{u\in U} A r_u$ and
$\Lambda = A + r = A + \sum _{u\in U}A r_u. $ This proves our claim
for $m=2$.

Assume now that the claim holds when $m \leq l $ ( where $ l \geq 2 $) and we show that the claim
also holds when $m= l+1.$ Let
$W$ denote the subalgebra generated by $A\cup B$ as algebras in $\Lambda $.
For $\bar \Lambda = \Lambda /r^l$, by Lemma \ref {3.2}, $\bar \Lambda = (A + r^l)/r^l \oplus r /r^l.$ It is clear
$r(\Lambda / r^l) = r/ r^l$. Indeed, obviously  $ r/ r^l\subseteq  r(\Lambda / r^l)$. Since
$(\Lambda / r ^l) /(r/r^l) \cong \Lambda / r$,  $r(\Lambda /r^l) \subseteq r/r^l$. Thus  $r(\Lambda / r^l) = r/ r^l$.
 Let
$\phi : \Lambda /r^2 \rightarrow (\Lambda /r^l)/(r^2 /r^l)$ be the canonical isomorphism, i.e.
$\phi (x + r^2) = (x+r^l) + (r^2/r^l)$ for any $x\in \Lambda .$
See
\begin {eqnarray*}
(r /r^l)/(r^2 /r^l)&=& \phi (r /r^2)\\
&=& \phi (\sum _{u\in U}( \Lambda r_u ) + r^2)\hbox { \ \ \ by assumption }\\
&=& (\sum _{u\in U} ( \Lambda  r_u  +r^l) + (r^2/r^l).
\end {eqnarray*}
Therefore,  $ \{  r_u  + r^l \mid u
\in U \} $ generates $(r/r^l)/ (r^2/ r^l) $ as $(\Lambda /r^l)/ (r/r^l)$-modules.
 By  induction assumption, we have
$\Lambda /r^l = (W +r^l)/r^l$.

Let $x \in \Lambda .$ There is $y \in W$ and $z \in r^l$ such that $x-y = z$.
Since $l \geq 2$, there exist
$\alpha _i \in r^{l-1},$ $\beta _i \in r$ for $i = 1, 2, \cdots ,n$ such that
$z = \sum \alpha _i \beta _i$.
Again using $\Lambda /r^l = (W +r^l)/r^l$, we have that
there are $a_i, b_i \in W, u_i, v_i \in r^l$ such that $\alpha _i = a_i + u_i$
and $\beta _i = b_i + v_i$,
so $a_i = \alpha _i -u_i \in r^{l-1}$ and $b_i = \beta _i - v_i\in r$ for any
$i = 1, 2, \cdots, n.$ By computation and $r^{l+1}=0$, we have $x-y \in W $
and $x \in W$. We complete the proof. $\Box$

Recall that $J$ is the ideal generated by all arrows in $D$ of $k(D, \Omega )$ and $\bar J$ is the
ideal $J/(\rho )$ of $k (D, \Omega , \rho)$.

\begin {Lemma}\label {3.6}{\ }\footnote  {The lemma was proved by Fang Li.}
If $J^t \subseteq (\rho )$ for some $t$, then $r (k(D , \Omega ,
\rho )) = \bar J$

\end {Lemma}

{\bf Proof.} Let $P= k (D, \Omega )$ and $Q = k(D, \Omega , \rho )$. Obviously  $Q/\bar J \cong
P/ J \cong \sum \{P_{ij}/ J_{ij} \mid i, j \in I\}$. It is clear that $P_{ij} = J_{ij}$ when
$i \not= j$ and $P_{ii}/ J_{ii} \cong \Omega _{ii}$. Thus $r (k(D , \Omega , \rho ))\subseteq \bar J $. Conversely,
since $J^t \subseteq (\rho )$ for some $t$, $\bar J$ is nilpotent and $\bar J
\subseteq r (k(D , \Omega , \rho )).$  $\Box$

\begin {Lemma}\label {3.7}

Let $\Lambda$ be an algebra.

(i) If $f$ is an algebra homomorphism from $k(D, \Omega )$ to $\Lambda $, then
$f$$\mid _{\Omega }$ is an algebra homomorphism
and $f(x_{ij}) = f(e_{ii})f(x_{ij})= f(x_{ij})f(e_{jj})$ for any arrow $x_{ij}$
from $i$ to $j$ and $i, j \in I$.

(ii) If $f$ is a map from $\Omega \cup D_1$ to $\Lambda $ and $f$$
\mid _ {\Omega }$ is an algebra homomorphism with $f(x_{ij}) =
f(e_{ii})f(x_{ij})= f(x_{ij})f(e_{jj})$ for any arrow $x_{ij}$
from $i$ to $j$ and $i, j \in I$, then there exists (unique)
algebra homomorphism $\bar f : k (D, \Omega ) \rightarrow \Lambda
$ such that $\bar f$$ \mid _{\Omega \oplus D_1} = f$.

\end {Lemma}

{\bf Proof. } (i) It is obvious.

(ii) Let $P$ denote the generalized path algebra $k (D, \Omega )$.
For any $i, j \in I$ and generalized path
$\alpha = a_{i_0}x_{i_0i_1} a_{i_1}x_{i_1i_2} \cdots a_{i_{n-1}}x_
{i_{n-1}i_n} a_{i_n} $
from $i_0=i$ to $i_n=j$, define $f_{ij}(\alpha ) =f(a_{i_0})f(x_{i_0i_1}) f
(a_{i_1})f(x_{i_1i_2})\cdots f(a_{i_{n-1}})f(x_{i_{n-1}i_n})f ( a_{i_n})$.
We get a $k$-linear map $f_{ij}$ from $P_{ij}$ to $\Lambda $.
Now we show
$$f_{is}(\alpha \beta ) = f_{ij}(\alpha ) f_{js} (\beta ) \ \ \ \ \ \ \ \ \ \ \
\ \ \ \ (*)$$

\noindent for two $\Omega $- paths \ \ \ \ \ \ \ \ \ \
$\alpha =a_{i_0} x_{i_0i_1}a_{i_1}x_{i_1i_2}a_{i_2}x_{i_2i_3} \cdots
x_{i_{n-1}i_{n}}a_{i_n}$
\ \ \ \ \ \ \ \ \ \ and \\
$\beta =b_{j_0}y_{j_0j_1} b_{j_1}y_{j_1j_2}b_{j_2}y_{j_2j_3} \cdots
y_{j_{m-1}j_{m}}b_{j_m}$ of $D$ with $i_n=j_0 = j, i_0 =i$ and $j_m = s$.
When $\alpha \beta \not=0,$ i.e. $a_{i_n} b_{j_0} \not=0,$ $\alpha \beta $ is
an $\Omega$-path. By definition, (*) holds.
When $\alpha \beta =0,$ i.e. $a_{i_n} b_{j_0} =0,$ $\alpha \beta $ is not an
$\Omega$-path. Obviously the left side of (*) =0.
\begin {eqnarray*}
\hbox { The right side of } (*) &=& f(a_{i_0})f(x_{i_0i_1}) f (a_{i_1})f
(x_{i_1i_2})\cdots f(a_{i_{n-1}})f(x_{i_{n-1}i_n})f( a_{i_n})\\
&{ \ }& f(b_{j_0})f(y_{j_0j_1})f( b_{j_1})f(y_{j_1j_2})f (b_{j_2})
f(y_{j_2j_3}) \cdots f(y_{j_{m-1}j_{m}})f(b_{j_m}) \\
&=& 0
\end {eqnarray*}

\noindent Consequently, (*) holds. For any $i,j \in I$, $f_{ij}$ naturally
becomes a $k$-linear map from $P_{ij}$ to $\Lambda $
with $f_{ij} (x_{is}y_{sj}) = f_{is}(x_{is})f_{sj}(y_{sj})$ and $f(x_{is}) = f
(e_{ii})f(x_{is})= f(x_{is})f(e_{ss})$
for any $x_{is} \in P_{is}$ and $y_{sj} \in P_{sj}$ and $i, s, j \in I$.
Let $\bar f =\oplus _{i, j \in I}f_{ij}.$ This $\bar f$ fulfills our requirement. $\Box$

Now we give our main theorem.

\begin {Theorem} \label {3.8} Algebra
$\Lambda $ can be lifted with nilpotent Jacobson radical  $r = r(\Lambda )$ and  has  gm unit
$\{e_{ii}'\}_I$ with each $\overline { e_{ii} '} $ in the center of $\bar \Lambda = \Lambda /r$ iff $\Lambda $ is
isomorphic to a generalized path algebra with weak relations.

\end {Theorem}

{\bf Proof.}
Assume that $\Lambda = A \oplus r$ with nilpotent Jacobson radical $r= r(\Lambda )$ and subalgebra $A$. By Lemma \ref {3.3},
 $e_{ii}' \in A$ for any $i\in I.$ Let
$e_{ii} = \overline  {e_{ii}'} = e_{ii}' +r$ in $\Lambda /r$ for any $i\in I.$
By Lemma \ref {3.1}, we have that $\pi \xi = id$, where
$\pi : \Lambda \rightarrow \Lambda /r(\Lambda )$ is the canonical homomorphism
and $\xi : \Lambda /r (\Lambda ) \rightarrow \Lambda $ is an algebra
homomorphism by defining
$\xi (x +r) = x_A$ for any $x = x_A +x_r\in \Lambda $ with $x_A \in A$ and $x_r \in r.$ Let
$\Omega _{ii}= e_{ii}(\Lambda /r)e_{ii}.$ Obviously $\{e_{ii}\}_I$ is gm unit of $\Omega $ and
$r(\Omega ) =0 $. For any
$i,j \in I$, let $ B_{ij} \subseteq e_{ii}' r e_{jj}'= r_{ij} $ such
that $\bar B_{ij} =: \{ \bar x = x + r ^2 \mid
x \in B_{ij } \} \subseteq r/r^2 $ is the $k$-basis of
$ \overline {e_{ii}'} (r/r^2) \overline {e_{jj}'} = e_{ii} (r/r^2) e_{jj}$.

We now construct a generalized path algebra $k (D , \Omega )$. Let $I$ be the
vertex set of $D$ and
$B_{ij}$  all of arrows from $i$ to $j$. Next we define an algebra
homomorphism $\varphi : k(D, \Omega )\rightarrow \Lambda $ by
$\varphi $$ \mid _\Omega = \xi $ and $\varphi (x) =x$ for
any arrow $x$ from $i$ to $j$.
Indeed, since $\xi (e_{ii}) = e_{ii}'$, we have $\varphi (x_{ij}) =x_{ij}$ and $\varphi (e_
{ii})\varphi (x_{ij})= \xi (e_{ii})\varphi (x_{ij})
= e_{ii}' x_{ij} = x_{ij}$, so $\varphi (x_{ij}) = \varphi (e_{ii}) \varphi (x_{ij})$ for any arrow
$x_{ij}$ from $i$ to $j$ and $i, j \in I$.
Similarly, $\varphi (x_{ij}) = \varphi (x_{ij})\varphi (e_{jj})$ for any arrow $x_{ij}$ from $i$ to
$j$ and $i, j \in I$.
By Lemma \ref {3.7}, $\varphi $ can become an algebra homomorphism from $k(D, \Omega )$ to $\Lambda $. Since $\bar
B_{ij}$ is a $k$-basis of $e_{ii}(r/r^2) e_{jj}$ for any $i, j \in I$ and
$r/r^2 = \sum _{i, j \in I} e_{ii}(r/r^2) e_{jj}$, $r/r^2$ is generated by
$\cup _{i, j \in I} \bar B_{ij}$ as $\Lambda /r$-modules. By Lemma \ref
{3.5}, $\Lambda$ is generated by $A \cup(\cup_{i, j\in I} B_{ij})$ as algebras.
This proves that $\varphi $ is surjective.

We now consider $N =: ker \varphi$.  Assume  $r^t =0$.
Since $\varphi (J) \subseteq r$, $\varphi (J^t) = 0$. Thus $J^t \subseteq
N$. For any $x\in ker \varphi $, obviously, there exist $a \in \Omega $ and
$\alpha \in J$
such that $ x= a +\alpha $. Thus $0 = \varphi (x) = \varphi (a) + \varphi (x) =
\xi (a) + \varphi (x) $. Considering  $\varphi (J) \subseteq r$ and
$\Lambda = A \oplus r$, we have $a =0$. $J^t \subseteq N \subseteq J$ has
been proved.

Conversely, assume that $\Lambda $ is a generalized path
algebra $k (D, \Omega , \rho )$ with weak relations. Let $P = k(D, \Omega),Q = k(D, \Omega , \rho )$ and $N = (\rho ).$
Since $P = \Omega \oplus J$ and $(\rho )\subseteq J$, by Lemma \ref{3.2},
we have that $Q = P/(\rho )= \Omega / (\rho ) \oplus J/ (\rho ).$  By Lemma \ref {3.6},
the Jacobson radical $r (Q) = \bar J$. Thus $Q$ can be lifted.  $r(Q)^t = \bar J ^t =0$ since
$J^t \subseteq
N.$
Since $\{e_{ii}\}_I$ is a complete set of pairwise orthogonal idempotents of $P$,
  $\{e_{ii} + N\}_I$ is a complete set of pairwise orthogonal idempotents of $Q$.
Obviously,  $\Omega \stackrel {\phi _1}{\cong } P /J \stackrel {\phi _2}{\cong }
Q/\bar J $  as algebras and $\phi _2\phi _1 (e_{ii}) = (e_{ii} +N) + \bar J $ for any $i \in I.$
 Since $e_{ii}$ is in the center of $\Omega $, $(e_{ii} +N) + \bar J$ is in center of $Q/ \bar J$ for any $i\in I.$
 $\Box$

\begin {Example}\label {3.8'}
Let $D$ be a directed graph with vertex set $I= {\bf N}$ of
natural numbers and only one arrow $x_{i, i+1}$ from $i$ to $i +1$
for any $i \in I.$ Let $\Omega _{ii}= M_{i}(k)$, the matrix
algebra of all $(i \times i )$- matrices  over $k$ for any $i\in
I$.  Set
$$ \rho = \{ x \mid  x = x_{i, i+1} x_{i+1, i+2}x_{i+2, i+3} \hbox {  is a path from } i  \hbox { to  }
i +3, i \in I\} \cup \{ x_{1, 2} \}.$$  Then $k(D, \Omega
)/(\rho)$ is a generalized path algebra with weak relations.
\end {Example}

\begin {Corollary} \label {3.9}  $\Lambda $ can be lifted with  nilpotent Jacobson radical
 and  with non-zero unity element iff $\Lambda $ isomorphic to a generalized path algebra  with one vertex
 and with weak relations

\end {Corollary}

{\bf Proof.} The sufficiency  follows from Theorem \ref {3.8} and its proof. We now show the necessity. Let $u$ be  the unity element of $\Lambda .$
Obviously, $\{u \}$ is a gm unit of $\Lambda $ and $\bar u$ is in the center of $\bar \Lambda = \Lambda /r(\Lambda ).$
By Theorem \ref {3.8} and its proof,  $\Lambda $ isomorphic to a generalized path algebra  $k (D, \Omega , \rho )$
with one vertex and with weak relations.
$\Box$

\begin {Lemma} \label {3.10} Let $\Lambda = A \oplus r$ with subalgebra $A$ and with nilpotent Jacobson radical
$r = r(\Lambda )$. If $\Lambda $ has the non-zero unity element $u$ and   $\{ \bar e_{ii}\}_I$ is a complete set
 of  pairwise orthogonal
idempotents of $ \bar \Lambda = \Lambda /r $, then $\{  (e_{ii})_A\}_I$ is a complete set of  pairwise orthogonal
idempotents of $  \Lambda $.
\end {Lemma}
{\bf Proof.} Let $\xi : \Lambda /r \rightarrow \Lambda $ by sending $x +r$ to $x_A$ for any $x\in \Lambda .$ Since
$\xi$ is an algebra homomorphism, we have that $\{  (e_{ii})_A\}_I$ is a  set of  pairwise orthogonal
idempotents. By Proposition \ref {1.3} (ii), $I$ is finite and $\bar u = \sum _{i\in I} \bar e_{ii}.$
By Lemma \ref {3.4}, $u \in A$.  Thus $ u = \sum _{i\in I}  (e_{ii})_A$  and
$\{  (e_{ii})_A\}_I$ is a complete set of  pairwise orthogonal
idempotents of $  \Lambda $. $\Box$

It is well known that, for any algebra $\Lambda $,  if $\Lambda / r(\Lambda )$ is a left (or right) artinian algebra
with non-zero unity element, then,
by  Wedderburn-Artin Theorem,
$\Lambda /r(\Lambda ) = B_1 \oplus B_2 \oplus \cdots \oplus B_n $ as algebras and
$B_i$ is a simple subalgebra of $\Lambda/ r(\Lambda ) $ for
any $i \in I = \{1, 2, \cdots , n\}.$ The number  $n$ is called the Wedderburn-Artin
number of $\Lambda $, written as  $n_{WA}(\Lambda ).$  If
$\Lambda / r(\Lambda )$ is not an artinian algebra with unity element, then  we write $n_{WA}(\Lambda ) = \infty.$

\begin {Corollary} \label {3.11}
(i) If $k(D, \Omega , \rho )$ is a generalized path algebra with weak relations, then $\mid $$D_0$$\mid \leq  n_{WA}
(k(D, \Omega , \rho ))$.

(ii) Let  $\Lambda $ can be lifted with  nilpotent Jacobson radical $r$  and
with non-zero unity element. If $\Lambda /r =  B_1 \oplus B_2 \oplus \cdots \oplus B_n $ as algebras and $B_i$ is a
non-zero
subalgebra
of $\Lambda /r(\Lambda )$  for   $i\in I  = \{1, 2, \cdots, n\},$
then  $\Lambda $ isomorphic to a generalized path algebra
$k(D, \Omega , \rho )$ with weak relations  and  $\Omega _{ii} =B_i$ for $i \in I =D_0.$

(iii) Let  $\Lambda $ can be lifted with  nilpotent Jacobson radical $r$  and
with non-zero unity element. If $\Lambda /r(\Lambda )$ is artinian, then  for any natural number
$m \leq  n_{WA}(\Lambda )$, $\Lambda $ isomorphic to
  a generalized path
algebra  $k(D, \Omega , \rho )$  with weak relations and
$\mid$$D_0$$\mid =m$.

\end {Corollary}
{\bf Proof.}
(i) Let $P = k (D, \Omega )$, $N = (\rho ) $ and $Q = P/ N.$ If $Q/r(Q)$ is artinian with unity element,  then,
by  Wedderburn-Artin Theorem,
$Q /r(Q ) = B_1 \oplus B_2 \oplus \cdots \oplus B_n $ as algebras and
$B_i$ is a simple subalgebra of $Q/ r(Q ) $ for
any $i \in  \{1, 2, \cdots , n\}.$
It is clear that $$ \oplus _{i\in I}\Omega _{ii}  \cong   B_1 \oplus B_2 \oplus \cdots \oplus B_n \ \ \ \
\hbox {as algebras. }$$
This implies that $$ \oplus _{i\in I}\Omega _{ii}  =   B_1' \oplus B_2' \oplus \cdots \oplus B_n' \ \ \ \
\hbox {as algebras },$$
where  $B_i'$ is a  simple subalgebra of $\Omega $ for $i=1, 2, \cdots , n.$
Considering  $B_1', B_2' , \cdots , B_n'$ are simple subalgebras, we have that each $\Omega _{ii}$ is a
sum of some of  $ \{ B_{1}', B_2' , \cdots , B_n' \}$. Thus $\mid$$I$$\mid  = \mid$$D_0$$\mid
 \leq n= n_{WA}(Q).$

If $Q/r(Q)$ is not an  artinian algebra with  the unity element, obviously $\mid$$D_0$$\mid
 \leq  n_{WA}(Q)$ since $ n_{WA}(Q) = \infty.$

(ii) Let $\Lambda = A \oplus r$ with subalgebra $A$  and $e_{ii}$  be the unity element of $B_i$ for any $i\in I.$ Obviously,
$\{ e_{ii}\}_I$ is a complete set of pairwise orthogonal central idempotents of $\Lambda /r.$ Let $e'_{ii}\in \Lambda $
such that $\overline {e_{ii}'} = e_{ii}$ for any $i\in I.$ By Lemma \ref {3.10}, $\{(e_{ii}')_A \}_I$ is a
complete set of pairwise orthogonal idempotents of $\Lambda $. By Theorem \ref {3.8} and its proof,
$\Lambda $ is isomorphic to   $k (D, \Omega , \rho )$ with weak relations  and  $\Omega _{ii} =B_i$ for $i \in I =D_0.$

(iii)  By  Wedderburn-Artin Theorem,
$\Lambda  /r(\Lambda ) = B_1 \oplus B_2 \oplus \cdots \oplus B_n $
as algebras and
$B_i$ is a simple subalgebra of $\Lambda / r(\Lambda  ) $ for
any $i \in  \{1, 2, \cdots , n\}$ with $n=n_{WA} (\Lambda ).$
Let $B_i' = B_i$ for $i=1, 2, \cdots , m-1$ and $B'_m = B_m + \cdots + B_n.$ Obviously,
 $\Lambda  /r(\Lambda ) = B_1' \oplus B_2' \oplus \cdots \oplus B_m' $
as algebras. By (ii),  $\Lambda $ is isomorphic to   $k (D, \Omega , \rho )$ with weak relations  and
$\mid$$D_0$$\mid =m$.
 $\Box$
\begin {Corollary} \label {3.12}

$\Lambda $ is isomorphic to a generalized path algebra with weak relations
when one of the following conditions holds:

(i) $\Lambda $ is a finite dimensional algebra with non-zero unity
element over a perfect field $k$ (e.g. the characteristic of $k$
is zero or $k$ is a finite field ).

(ii) $\Lambda $ is a finite-dimensional separable algebra with non-zero unity element.

(iii) $\Lambda $ is an algebra over a field $k$ with non-zero unity element and nilpotent
Jacobson radical, and $sup \{ n \mid H_k^n (\Lambda , M) \not=0$ for some
$\Lambda $-bimodule $M \}\leq 1$ (see \cite [Definition 11.4] {Pi82}).

\end {Corollary}

{\bf Proof.} It follows from the famous Wedderburn-Malcev Theorem (see \cite
[Theorem 11.6 and Corollary 11.6] {Pi82}) that $\Lambda $ can be lifted.
We complete the proof by
 Corollary  \ref {3.9}. $\Box$

\begin {Corollary} \label {3.13}

Let $k$ be a perfect field.

(i)  $\Lambda $ is a finite dimensional algebra with non-zero unity element iff
$\Lambda $ is isomorphic to a generalized path
algebra $k(D, \Omega , \rho )$ of finite directed graph  with weak relations and with $dim \ \Omega < \infty $.

 (ii) If $\Lambda $ is a finite dimensional algebra with non-zero unity element over field $k$, then
$\Lambda $ is isomorphic to a generalized path
algebra  $k (D, \Omega , \rho )$  of finite directed graph with weak relations and  $\Omega _{ii}= B_i$ for any $i\in I= \{1, 2, \cdots , n\}$.
Here

 $\Lambda /r = B_1 \oplus B_2 \oplus \cdots \oplus B_n $ as algebras and $B_i$ is a simple subalgebra of $\Lambda /r$ for
any $i \in I.$

(iii) If $\Lambda $ is a finite dimensional algebra with non-zero unity element over field $k$, then  for any natural number
$m \leq  n_{WA}(\Lambda )$, there exists  a generalized path
algebra  $k(D, \Omega , \rho )$  with weak relations and $\mid$$D_0$$ \mid =m$.

\end {Corollary}

{\bf Proof.} (i) $\Lambda $ is a finite dimensional  algebra with non-zero unity element over field $k$, then
$\Lambda $ is isomorphic to a generalized path
algebra of finite directed graph  with weak relations and $dim {\ } \Omega < \infty $ by corollary \ref {3.12} and the proof of Theorem
\ref {3.8}.
Conversely, assume
$\Lambda = k(D, \Omega , \rho )$ is  a generalized path
algebra of finite directed graph  with weak relations. Let  $P = k(D, \Omega )$,  $Q = k(D, \Omega , \rho )$
and $N = (\rho )$.
For any $i, j\in I$, $Q_{ij}$ is spanned by $ \{  [\alpha ] + N \mid  \alpha  $ is a generalized path from $i$ to $j$
with $l(\alpha ) \leq  t  \}$ since $J^t \subseteq (\rho )$. However,  $ \{  [\alpha ]
 \mid  \alpha  $ is a generalized path from $i$ to $j$
with $l(\alpha ) \leq  t  \}$  is spanned by finite elements  since  $\Omega $ is
finite dimensional.
Consequently, $Q$ is finite dimensional.

(ii) By \cite [Corollary 11.6] {Pi82}, $\Lambda $ can be lifted.  Obviously the Jacobson
radical $r$ is nilpotent. By  Wedderburn-Artin Theorem,
$\Lambda /r = B_1 \oplus B_2 \oplus \cdots \oplus B_n $ as algebras and $B_i$ is a simple subalgebra of $\Lambda $ for
any $i \in I = \{1, 2, \cdots , n\}.$ Using Corollary \ref {3.11}(ii), we complete the proof.

(iii) It follows from Corollary \ref {3.11}(iii) and  \cite [Corollary 11.6] {Pi82}. $\Box$

\begin {Example} \label {3.13'} Let $k$ be the complex field or
real field and $\Lambda$ the matrix algebra $M_n (k)$ of all $(n
\times n )$- matrices  over $k$. Then it follows from Corollary
\ref {3.13} that $\Lambda $ is isomorphic to a generalized path
algebra $k(D, \Omega , \rho )$ of finite directed graph  with weak
relations and with $dim \ \Omega < \infty $.
\end {Example}

An algebra $\Lambda $ over field $k$ is called a generalized elementary algebra
if $\Lambda / r (\Lambda )\cong \oplus _{i\in I} B_{ii}$ as algebras with $B_
{ii} = k$ for any $i \in I$. A finite dimensional generalized
elementary algebra with unity element is called an elementary algebra.

\begin {Corollary} \label {3.14}  $\Lambda $ is a generalized elementary algebra  which can be lifted
 with nilpotent Jacobson radical
$r = r(\Lambda )$ and  has  a complete   set of  pairwise orthogonal
idempotents
iff $\Lambda $ is isomorphic to a path algebra with relations.

\end {Corollary}

{\bf Proof.} The sufficiency follows from Theorem \ref {3.8}.
We now show the necessity. Assume that $\Lambda = A \oplus r $ and
$\Lambda /r = \oplus _{i\in I} k \bar e_{ii}$ as algebras, where $A$ is a subalgebra of $\Lambda $ and $r$ is the Jacobson
radical of $\Lambda .$  Obviously, $\{\bar e_{ii}\}_I$ is a
complete set of pairwise orthogonal central idempotents of $\bar \Lambda = \Lambda /r.$
Let $\xi : \Lambda /r \rightarrow \Lambda $ by sending $x +r$ to $x_A$ for any $x\in \Lambda .$ Since
$\xi$ is an algebra homomorphism by Lemma \ref {3.1}, we have that $\{  (e_{ii})_A\}_I$ is a  set of  pairwise orthogonal
idempotents. However, $\Lambda = (\sum _{i\in I} k(e_{ii})_A) + r$. For any $x \in (\sum _{i\in I} k(e_{ii})_A) \cap r$,
there exist $\alpha _i \in k$ such that $x = \sum _{i\in I} \alpha _i(e_{ii})_A $.
 Since  $0 = \bar x = \sum _{i\in I} \alpha _i \overline { (e_{ii})_A } $, we have $\alpha _i = 0$ for any $i\in I$. This implies
$x=0$ and $\Lambda = (\sum _{i\in I} k(e_{ii})_A) \oplus r$. Since $(\sum _{i\in I} k(e_{ii})_A) \subseteq A$,
$\sum _{i\in I} k(e_{ii})_A=A$.

Let $\{  e_{ii}'\}_I$ be a complete  set of  pairwise orthogonal
idempotents of $\Lambda $. By Lemma \ref {3.3}, $\{  e_{ii}'\}_I\subseteq A =\sum _{i\in I} k(e_{ii})_A.$
Since $\{  e_{ii}'\}_I$ is a complete set then so is  $\{ ( e_{ii})_A \}_I$.
By Theorem \ref {3.8}, $\Lambda $ is isomorphic to a path algebra with weak relations.

It remains to show $ ker \varphi \subseteq J^2, $ where $\varphi $ is the same as in the proof of Theorem \ref {3.8}.
For any $x\in ker \varphi $, obviously, there exist $y \in J $, $y \not\in J^2$
and $z \in J^2$ such that $ x= y + z$. Thus
$0 = \varphi (x) = \varphi (y) + \varphi (z) $ and $\varphi (z) \in r^2$. Thus  $\varphi (y) \in r^2.$ Since
$y \in J$ and $y \not\in J^2,$  there are  mutually different  arrows $x_1, x_2, \cdots , x_n $ such that
$y = \sum _{p =1} ^n  \alpha _p x_p$ with
$\alpha _p \in k$ for $p = 1, 2, \cdots , n.$
Notice $x_1, x_2, \cdots , x_n \in  \cup _{i, j \in I}B_{ij}$, where $B_{ij}$  is
 the same as in the proof of Theorem \ref {3.8}.
See that $0 = \overline  { \varphi (y) } =
 \sum _{p =1} ^n  \alpha _p \bar  x_p$  in $r/r^2.$ However, $ \{ \bar x_1,  \bar x_2, \cdots , \bar  x_n \}$  is
independent, so $\alpha _p =0 $ for $p = 1, 2, \cdots , n.$ This implies $y=0.$
 Consequently, $ker \varphi \subseteq J^2.$ $\Box$

There exist generalized elementary algebras whose Jacobson radicals are not nilpotent.

\begin {Example} \label {3.15}

Let $D$ be a directed graph with vertex set $I= {\bf N}$ of natural numbers and
only one arrow from $i$ to $i +1$ for any $i \in I.$ Path algebra $kD$ is
an elementary algebra since its Jacobson radical $r(kD)$ is $J$. However, $r(kD)$ is not nilpotent.

\end {Example}

It immediately follows from Corollary \ref {3.14} that
\begin {Corollary} \label {3.16}  $\Lambda $ is an  elementary algebra  which can be lifted
 iff $\Lambda $ is isomorphic to a path algebra of finite directed graph with relations.

\end {Corollary}

{\bf Remark:} In the above corollary, we require the condition that $\Lambda $ can be lifted, but
this was not mentioned explicitly  in
\cite [Theorem 1.9] {ARS95}.
Assume that $\Lambda /r = \oplus _{i = 1, 2, \cdots , n} k \bar e_{ii}$ as algebras. It is clear that
 there exists a complete set $\{ e_{ii}' \mid i = 1, 2, \cdots , m \}$  of pairwise
orthogonal primitive idempotemts of $\Lambda $.  In the proof of  \cite [Theorem 1.9] {ARS95},  the condition $m=n$
was used without
proof. However, this condition implies that $\Lambda $ can be lifted.  Indeed, since $e_{ii}'$ is non-zero idempotent,
$e_{ii}' \not\in r$ for any $i= 1, 2, \cdots, n$. Thus $\{ \bar {e_{ii}'}  \mid i = 1, 2, \cdots , n \}$  is
linear independent in $\bar \Lambda = \Lambda /r.$  Consequently, $\Lambda /r =
\oplus _{i = 1, 2, \cdots , n} k \bar e_{ii} = \oplus _{i = 1, 2, \cdots , n} k \bar {e_{ii}'}$. It is easy to check
  $\Lambda  =
(\oplus _{i = 1, 2, \cdots , n} k  e_{ii}') \oplus r$ and $(\oplus _{i = 1, 2, \cdots , n} k  e_{ii}')$ is a subalgebra
 of $\Lambda.$ That is, $\Lambda $ can be lifted. $\Box$

Finally we give  the gradations of the gm algebras and the generalized path
algebras.

\begin {Proposition} \label {3.17} (see \cite [Proposition 2.1] {ZZ03a})
Let  $A=\sum \{A_{ij} \mid i, j\in I \}$ be  a gm algebra and  $G$
 an abelian group. If there exists a bijective map $\phi : I \rightarrow  G,$
 then $A$ is an algebra graded by $G$ with
$A_g = \sum _{\phi (i) = \phi (j) +g} A_{ij}$ for any $g\in G$. In this case, the gradation is called a generalized matrix gradation,
or gm gradation in short.
\end {Proposition}
{\bf Proof.} For any $g, h \in G$, see that
\begin{eqnarray*}
A_g A_h &=& (\sum _{\phi (i) = \phi (j) +g} A_{ij})( \sum _{\phi (s) =\phi ( t)+h}
A_{st})\\
 &\subseteq & \sum _{\phi (i) = \phi(t) + h +g}A_{i,\phi ^{-1}(\phi (t) +h))  }A_{ \phi ^{-1}(\phi (t) +h), t}\\
 &\subseteq & A_{g +h}. \ \ \ \ \ \
 \end{eqnarray*}
Thus $A=\sum\{A_{ij} \mid i, j\in I \}= \sum _{g\in G} A_g$ is a
$G$-grading algebra. $\Box$

\begin {Proposition} \label {3.18}

(i) Let $Q = k (D, \Omega , \rho )$ be a
generalized path algebra with weak relations.
If $D_0$ is finite, then $Q$ has  a gm gradation  by ${\bf Z}_m$ when $m \leq D_0$.

(ii) Assume that  $\Lambda $ can be lifted with  nilpotent Jacobson radical $r$  and
with non-zero unity element. If $\Lambda /r(\Lambda )$ is artinian, then  for any natural number
$0\not= m \leq  n_{WA}(\Lambda )$, $\Lambda $  has  a gm gradation  by ${\bf Z}_m$.

(iii) If $\Lambda $ is a finite dimensional algebra with non-zero unity element over perfect field $k$, then  for any natural number
$m \leq  n_{WA}(\Lambda )$, $\Lambda $  has  a gm gradation  by ${\bf Z}_m$.

\end {Proposition}

{\bf Proof.}
(i) Assume $D_{0}= \{1, 2, \cdots , n\}$.  Let $e_{ii}' = e_{ii}$  for $i= 1, 2, \cdots m-1$, $e_{mm}' = e_{mm} +
\cdots  + e_{nn}$. It is clear that $\{e_{ii}'\}$ is a complete set of pairwise orthogonal idempotents of $Q$
with  $\bar e_{ii}'$ in the center of $Q / r(Q)$ since $\bar e_{jj}$ is in  the center of $Q / r(Q)$  for any $i
= 1, 2, \cdots m$ and $ j = 1, 2, \cdots , n $. By Theorem \ref {3.8}, $Q$ can be lifted. It follows from
Theorem \ref {3.8} that $Q$ is isomorphic to a generalized path algebra with weak relations and with $m$ vertexes.
 By Proposition \ref {3.17}, $Q$ has  a gm gradation
by ${\bf Z}_m$.

(ii) It follows from Proposition \ref {3.17} and Corollary \ref {3.11} (iii).

(iii) It follows Corollary \ref {3.13} and Proposition \ref {3.17}. $\Box$

\vskip 1.0cm

\noindent {\bf Acknowledgement }: The work was supported by
Australian Research Council. The first author thanks the
Department of Mathematics, University of Queensland for
hospitality. He also  wishes to thank F. Li for assistance.

\end {document}